\documentclass[12pt]{article}
\usepackage{amsmath, amssymb, amsfonts}
\usepackage{cite}

\newtheorem{theorem}{Theorem}
\newtheorem{lemma}{Lemma}

\newcommand{\proof}{\textit{Proof.}}
\newcommand{\qed}{$\Box$}
\newcommand{\op}{\operatorname}
\newcommand{\abs}[1]{\left\lvert #1 \right\rvert}
\newcommand{\icard}[2]{\abs{#1 \cap #2}} 
\newcommand{\floor}[1]{\left\lfloor #1 \right\rfloor}
\newcommand{\irange}[2]{\left\{#1,\ldots,#2\right\}} 

\begin{document}

\begin{center}

{\LARGE
Large cycles in random generalized Johnson graphs
}

{\Large
V.S.~Kozhevnikov\footnote[1]{Moscow Institute of Physics and Technology (National Research University), Department of Discrete Mathematics, Dolgoprudny, Moscow Region, Russian Federation}, A.M.~Raigorodskii\footnotemark[1]\footnote[2]{Lomonosov Moscow State University, Faculty of Mechanics and Mathematics, Moscow, Russian Federation}\footnote[3]{Adyghe State University, Caucasus mathematical center, Maykop, Republic of Adygea, Russian Federation}\footnote[4]{Institute of Mathematics and Computer Science, Buryat State University, Ulan-Ude, Buryat Republic, Russian Federation}, M.E.~Zhukovskii\footnotemark[1]\footnotemark[3]\footnote[5]{The Russian Presidential Academy of National Economy and Public Administration, Moscow, Russian Federation}\footnote[6]{Moscow Center for Fundamental and Applied Mathematics, Moscow, Russian Federation}
}

\vspace{0.5cm}

Abstract\\

\end{center}

This paper studies thresholds in random generalized Johnson graphs for containing large cycles, i.e. cycles of variable length growing with the size of the graph. Thresholds are obtained for different growth rates.

\vspace{0.3cm}

{\bf Keywords:} random graphs, Johnson graph, Kneser graph, large cycles, threshold

\vspace{0.5cm}

\section{Introduction and new results}
\label{intro}

A simple graph $G(n,r,s) = (V, E)$ is called a \textit{generalized Johnson graph} if
\begin{enumerate}
\item $0\le s < r < n$,
\item $V = \binom{[n]}{r}$ is the set of all $r$-subsets of the set $[n] := \irange{1}{n}$,
\item $\forall x,y \in V : \left\{x,y\right\} \in E \Leftrightarrow \icard{x}{y} = s$.
\end{enumerate}
Note that the special case $G(n,r,r-1)$ is known as the Johnson graph and $G(n,r,0)$ is known as the Kneser graph.

Throughout the article it is assumed that $r$ and $s$ are constant and $n$ approaches $+\infty$. The total number of vertices in this graph is denoted by $N$:
\[N = \abs{V} = \binom{n}{r} \sim \frac{n^r}{r!}.\]
From the symmetry of the definition it is evident that $G(n,r,s)$ is a regular graph. Let $N_1$ denote the degree of its vertex: \[N_1 = \binom{r}{s} \binom{n-r}{r-s} \sim \binom{r}{s} \frac{n^{r-s}}{(r-s)!}.\]

A \textit{random generalized Johnson graph} $G_p(n,r,s)$ \ is a binomial random subgraph of $G(n,r,s)$. It is obtained from $G(n,r,s)$ by independent removal of each edge with probability $(1 - p)$. This is a generalization of the classical Erd\H{o}s-R\'{e}nyi random graph $G(n,p)$ \cite{Erdos1959, Gilbert1959}, corresponding to the case $r=1, s=0$, in which $G(n,1,0) \cong K_n$. Several recent works \cite{Bollobas2016, Kupavskii2016, Kupavskii2018, Zhukovskii2010, Zhukovskii2012_zero, Zhukovskii2012_sub, Burkin2016, Pyaderkin2016, Burkin2018, Ogarok2020, Raigorodskii2021} are devoted to the study of asymptotic properties of $G_p(n,r,s)$.

Let $\mathcal{A}$ be some graph property. As usual, we write $G \vDash \mathcal{A}$ when $\mathcal{A}$ holds for $G$. The function $\hat{p}=\hat{p}(n)$ is called a \textit{threshold} for the \textit{increasing} (see \cite{Janson2000}) property $\mathcal{A}$ if
\[\mathbb{P}\left(G_p(n,r,s)\vDash\mathcal{A}\right) \to
\begin{cases}
0, p = o(\hat{p}), \\
1, p = \omega(\hat{p}),
\end{cases}\]
where $p = \omega(\hat{p})$ means $\hat{p} = o(p)$. A threshold $\hat{p}$ is called \textit{sharp} if $\forall \varepsilon >0$:
\[\mathbb{P}\left(G_p(n,r,s)\vDash\mathcal{A}\right) \to
\begin{cases}
	0, p < (1-\varepsilon) \hat{p}, \\
	1, p > (1+\varepsilon) \hat{p}.
\end{cases}\]
The problem of finding thresholds for various properties and, more generally, the asymptotic behavior of probabilities of random graphs' properties has been widely studied \cite{Schurger1976, Komlos1983, Bollobas1987, Friedgut2006, Erdos1960, Bollobas1981, Karoski1983, Rucinski1986, Rucinski1988}.

Particularly well understood is the property of subgraph containment in $G(n,p)=G_p(n,1,0)$. Thresholds for containing a copy of an arbitrary fixed graph are found in classical works \cite{Erdos1960, Bollobas1981, Rucinski1986}.
Asymptotic Poisson distribution for the number of copies of a fixed strongly balanced graph inside the threshold is established in \cite{Bollobas1981}. The number of copies of a balanced (but not strongly balanced) graph is also known to be convergent to some distribution, but not necessarily Poisson \cite{Bollobas1989}. 
Asymptotic normality of the number of copies of a fixed graph for large enough $p$ is determined in \cite{Karoski1983, Rucinski1986, Rucinski1988}. Most of these results are generalized by Burkin for $G_p(n,r,s)$ \cite{Burkin2016}.

In this work we go further and consider subgraphs of variable size (i.e. depending on $n$) in $G_p(n,r,s)$. Specifically, we focus on thresholds for the property of variable-length cycle containment in $G_p(n,r,s)$. In \cite{Burkin2016} Burkin found (under certain assumptions) the threshold for the property $\mathcal{A}_H$ of containment in $G_p(n,r,s)$ of a subgraph isomorphic to a fixed graph $H$. Applied to a simple cycle $C_t$ of a fixed length $t$, his theorem can be stated as follows.

\begin{theorem}[Burkin, 2016, \cite{Burkin2016}]
Let $0 \le s < r$ and $t$ be fixed integers. Then the threshold for containment of $C_t$ in $G_p(n,r,s)$ is
\begin{equation}
\hat{p} = n^{-(r - s) - s / t}.
\label{eq_fixed_thresh}
\end{equation}
\label{th_fixed}
\end{theorem}

In this paper we generalize this result to growing cycles, i.e. $t\to+\infty$ as $n\to+\infty$. Moreover, if the growth rate of $t$ is fast enough, the threshold is shown to be sharp. These results are summarized in the following theorem.

\begin{theorem}
Let $0 \le s < r$ be fixed integers. If $t=t(n)$ satisfies the condition
\begin{equation}
t = o\left(\sqrt{N_1}\right),
\label{eq_t_o_n_1}
\end{equation}
then the threshold for containment of $C_t$ in $G_p(n,r,s)$ is
\begin{equation}
\hat{p} = \frac{n^{-s/t}}{N_1}.
\label{eq_sharp_thresh}
\end{equation}
Moreover, if $s=0$ and $t\to+\infty$ or $s$ is arbitrary and $t=\omega(\ln n)$, then the threshold \eqref{eq_sharp_thresh} is sharp.
\label{th_long_thresh}
\end{theorem}

Note that if $t = \op{const}$, then there is no sharp threshold, which follows from the fact that in this case, for $p = c \cdot n^{-(r-s) - s/t}$, where $c>0$ is any constant, the number of cycles has asymptotically Poisson distribution \cite{Burkin2016}.

The proof of Theorem~\ref{th_long_thresh} presented in Section~\ref{proofs} uses the typical first and second moment methods. In contrast to $t=\op{const}$, the case $t\to+\infty$ requires more careful bounds of the expectation and the variance of the number of cycles.

\section{Proofs}
\label{proofs}

The proof of Theorem~\ref{th_fixed} provided in \cite{Burkin2016} by Burkin is also based on the methods of the first and the second moments. Let $X$ and $c_t$ be the numbers of copies of $C_t$ in $G_p(n, r, s)$ and in $G(n, r, s)$ respectively. In order to estimate the expectation $\op{E} X = c_t p^t$, Burkin finds the asymptotics of $c_t$. For this purpose, he shows that $c_t$ asymptotically equals the number of cycles entirely contained in the neighborhood of a single vertex, which is easy to compute explicitly. However, this is not generally true if $t\to+\infty$, in which case the estimation of $c_t$ should be more subtle. This problem is resolved in Section~\ref{proofs_c_t}.

Estimation of the variance $\op{Var} X$ for $t=\op{const}$ is almost the same as in the Erd\H{o}s-R\'{e}nyi model for any balanced subgraph. If, however, $t\to+\infty$, then this problem becomes non-trivial. In Section~\ref{proofs_var}, a sufficient upper bound for $\op{Var} X$ is obtained using the specificities of $C_t$ topology.

Finally, using the bounds for $\op{E} X$ and $\op{Var} X$, Theorem \ref{th_long_thresh} is proved in Section~\ref{proofs_proof}.

\subsection{The number of cycles}
\label{proofs_c_t}

Counting simple cycles in $G(n,r,s)$ can be reduced to counting simple paths due to the following lemma.

\begin{lemma}
	Let $\{x, y\} \in E$ be an arbitrary edge in $G(n,r,s)$. Then the number $p_t$ of simple paths on $t$ vertices, whose first vertex is $x$ and last vertex is $y$, does not depend on $x$ and $y$ and
	\begin{equation}
		c_t = \frac{1}{2t}N \cdot N_1 \cdot p_t.
		\label{eq_c_t_in_terms_of_p_t}
	\end{equation}
	\label{lemma_c_t_in_terms_of_p_t}
\end{lemma}
\proof{}
The fact that $p_t$ does not depend on $x$ and $y$ is evident from the symmetry of $G(n,r,s)$ definition.

$C_t$ can be embedded in $G(n,r,s)$ in the following way. First, choose an arbitrary vertex $y \in V$ in $N$ ways. Next, choose its neighbor $x$ in $N_1$ ways. Then, choose a simple path on $t$ vertices with $x$ as its first vertex and $y$ as its last vertex in $p_t$ ways. Finally, since the number of automorphisms of $C_t$ is $2t$, the equality~\eqref{eq_c_t_in_terms_of_p_t} is obtained.
\qed{}

\bigskip
We now switch to estimating $p_t$. For $y\in V$, consider sets
\[V_j(y) = \left\{x \in V \mid \icard{x}{y} = j\right\}, j \in \irange{0,1}{r},\]
which form a partition of $V$. Let $ i, j \in \irange{0,1}{r}$, $y \in V$, and $x \in V_i(y)$. Following Burkin \cite{Burkin2016}, we define
\[A_i^j = \abs{V_j(y) \cap V_s(x)},\]
which does not depend on the choice of $x$ and $y$ due to the symmetry of $G(n,r,s)$ definition.

Note that if $i = r$, then $\icard{x}{y} = r$, which means that $x = y$. Therefore, $A_{r} ^ {j} = {{\delta}_{sj}} \cdot {{N}_{1}} $, where ${{\delta}_{sj}}$ is the Kronecker delta. If $ j = r $, then $\forall z\in V_j(y)$ we have $\icard{z}{y} = r $, which means that $z=y$ and, therefore, that $ A_ {i} ^ {r} = {{\delta}_{is}}$. For all $i,j \le r$, the following formula holds:
\begin{equation}
	\begin{gathered}
		A_{i}^{j}=\sum\limits_{m=0}^{s} \binom{i}{m} \binom{r-i}{s-m} \binom{r-i}{j-m} \binom{n-2r+i}{r-s-j+m} = \\ = \sum\limits_{m={{m}_{min}}}^{{{m}_{max}}} \binom{i}{m} \binom{r-i}{s-m} \binom{r-i}{j-m} \binom{n-2r+i}{r-s-j+m}
	\end{gathered}
	\label{eq_a_i_j_sum}
\end{equation}
assuming that $\binom{0}{0} \equiv 1$, $\binom{k}{k} \equiv \binom{k}{<0} \equiv \binom{<0}{k} \equiv 0$ for any $k \in \mathbb{Z}$. ${m}_{min}$ and ${m}_{max}$ are respectively the minimum and the maximum values of $m$ for which the corresponding term in the sum is nonzero.

It is clear that the $m$th term in the sum~\eqref{eq_a_i_j_sum} is nonzero iff
\[\max \{0,\max \{i,j\}-(r-s),i+j-r\}\le m\le \min \{s,i,j,n-2r+i+j-(r-s)\}.\]
Therefore, for sufficiently large $n$:
\begin{equation}
	\left\{ \begin{aligned}
		& {{m}_{min}}=\max \{0,\ \max \{i,j\}-(r-s),\ i+j-r\}, \\ 
		& {{m}_{max}}=\min \{s,i,j\}. \\ 
	\end{aligned} \right.
	\label{eq_m_min_max}
\end{equation}
The necessary and sufficient condition on $i$ and $j$ under which $A_{i}^{j} \ne 0$ (for sufficiently large $n$) is ${{m}_{min}}\le {{m}_{max}}$, or, equivalently,
\begin{equation}
	\left\{ \begin{aligned}
		& |i-j|\le r-s, \\ 
		& i+j\le r+s. \\ 
	\end{aligned} \right.
	\label{eq_when_a_i_j_nonzero}
\end{equation}
Note that for $i, j$ satisfying \eqref{eq_when_a_i_j_nonzero} the value ${m}_{max}$ defines the asymptotic order of $A_{i}^{j}$, namely,
\begin{equation}
	A_{i}^{j}\sim\operatorname{const}\cdot {{n}^{r-s-j+{{m}_{max}}}}=\operatorname{const}\cdot {{n}^{r-s-j+\min \{s,i,j\}}}.
	\label{eq_a_i_j_asymp}
\end{equation}
In particular,
\begin{equation}
	A_s^s\sim\frac{n^{r-s}}{(r-s)!}\sim\frac{N_1}{\binom{r}{s}}
	\label{eq_a_s_s_asymp}.
\end{equation}

\begin{lemma}
	If $t = o \left(\sqrt{N_1}\right)$, then
	\begin{equation}
		\begin{cases}
			p_t = O \left(N_1^{t-2}\right), \\
			p_t = \Omega\left(\left(N_1 \middle/ \binom{r}{s}\right)^{t-2}\right), \\
			p_t = \Omega\left(n^{-s} N_1^{t-2}\right).
		\end{cases}
		\label{eq_p_t_upper_and_two_lower_bounds}
	\end{equation}
	\label{lemma_p_t_upper_and_two_lower_bounds}
\end{lemma}
\proof{}
For $t=\op{const}$, the asymptotics of $p_t$ is derived in \cite{Burkin2016}:
\[p_t \sim \left(A_s^s\right)^{t-2}.\]
Then, \eqref{eq_p_t_upper_and_two_lower_bounds} follows from \eqref{eq_a_s_s_asymp}.
Therefore, it is enough to assume that $t$ is greater than any fixed constant for sufficiently large $n$.

To obtain the bounds in \eqref{eq_p_t_upper_and_two_lower_bounds}, let's fix some adjacent vertices $x$ and $y$ and estimate the number of ways to choose a simple path $x_1, \ldots , x_t$, where $x_1 = x$ and $x_t = y$. We will choose the vertices $x_2, \ldots, x_{t-1}$ sequentially so that, for $2\le j\le t$, $\{x_{j-1}, x_j\}\in E$, or, equivalently, $x_j \in V_s(x_{j-1}),$ occasionally imposing additional restrictions on the choice of $x_j$.

The upper bound in \eqref{eq_p_t_upper_and_two_lower_bounds} follows from the obvious fact that
\[p_t \le N_1^{t-2}.\]

The first lower bound in \eqref{eq_p_t_upper_and_two_lower_bounds} follows from \eqref{eq_a_s_s_asymp} and the fact that
\[p_t \ge (A_s^s - t)^{t-2},\]
which can be obtained by choosing only $x_2, \ldots , x_{t-1} \in V_s(y)$.
\bigskip

It remains to prove the last equality in \eqref{eq_p_t_upper_and_two_lower_bounds}. Given no other restrictions, the vertices $x_2, \ldots, x_{t - s - 1}$ can be chosen in no less than \[(N_1 - t)^{t - s - 2}\sim N_1^{t - s - 2}\] ways (assuming that for each such choice there exists at least one relevant choice of the rest of the vertices, which is shown below). Having fixed $x_2, \ldots, x_{t - s - 1}$, let's restrict the choice of $x_{t-s}, \ldots, x_{t-1}$ as follows:
\[x_{t-s-1+j} \in V_{\min\{i+j,s\}}(y),\]
where $1\le j\le s$, $i = \icard{x_{t - s - 1}}{y}$, i.e. $x_{t - s - 1} \in V_i(y)$. Then $x_{t-s}, \ldots, x_{t-1}$ can be chosen in no less than
\begin{equation}
	\label{_eq_the_rest_bound}
	\left(A_i^{\min\{i+1,s\}} - t\right) \cdot \left(A_{\min\{i+1,s\}}^{\min\{i+2,s\}} - t\right) \cdot\ldots\cdot \left(A_{\min\{i+s-1,s\}}^{\min\{i+s,s\}} - t\right) \tag{$*$}
\end{equation}
ways. If $r-s \ge 2$, then, applying \eqref{eq_a_i_j_asymp}, for $j \in \irange{0,1}{s-1}$:
\[A_{j}^{j+1} = \Theta\left(n^{r-s-1}\right) = \Omega\left(n^{(r-s)/2}\right) = \omega(t)\]
and
\[A_s^s = \Theta\left(n^{r-s}\right) = \omega(t).\]
Thus, for $r-s \ge 2$, the expression in \eqref{_eq_the_rest_bound} equals $\Omega \left(n^{(r-s-1)\cdot s}\right)$. Therefore,
\[p_t = \Omega\left(N_1^{t-s-2} \cdot n^{(r-s-1)\cdot s}\right) = \Omega\left(n^{-s} N_1^{t-2}\right).\]

If, however, $r-s=1$, then $A_{j}^{j+1} = \op{const}$, which makes \eqref{_eq_the_rest_bound} meaningless for large $t$.
Let's adopt another approach in this case choosing the vertices only from the following vertex sets:
\[\begin{gathered}
	x_{2} \in V_{s-1}(y), x_{3} \in V_{s-2}(y), \ldots , x_{s+1} \in V_{0}(y), \\
	x_{s+2}, \ldots, x_{t-s-1} \in V_{0}(y), \\
	x_{t-s} \in V_{1}(y), x_{t-s+1} \in V_{2}(y), \ldots , x_{t-1} \in V_{s}(y),
\end{gathered}\]
which can be done in no less than
\[\begin{gathered}
	A_s^{s-1} \cdot A_{s-1}^{s-2} \cdot\ldots\cdot A_{1}^{0} \cdot \left(A_0^0 - t\right)^{t-2s-2} \cdot \left(A_{0}^{1} - 1\right) \cdot \left(A_{1}^{2} - 1\right) \cdot\ldots\cdot \left(A_{s-1}^{s} - 1\right) \ge \\
	\ge \Theta(N_1^s) \cdot N_1^{t-2s-2} = \Theta\left(N_1^{t-2}N_1^{-s}\right) = \Theta\left(n^{-s} N_1^{t-2}\right)
\end{gathered}\]
ways, where the inequality follows from \eqref{eq_a_i_j_asymp}.

Thus, in any case,
\[p_t = \Omega\left(n^{-s} N_1^{t-2}\right).\]
The last bound in \eqref{eq_p_t_upper_and_two_lower_bounds} is proved.
\qed{}

\begin{lemma}
	If $t = o \left(\sqrt{N_1}\right)$, then
	\begin{equation}
		\begin{cases}
			c_t = O \left(\frac{n^s}{2t} N_1^{t}\right), \\
			c_t = \Omega\left(\frac{n^s}{2t} \left(N_1 \middle/ \binom{r}{s}\right)^{t}\right), \\
			c_t = \Omega\left(\frac{1}{2t} N_1^{t}\right). \\
		\end{cases}
		\label{eq_c_t_upper_and_two_lower_bounds}
	\end{equation}
	\label{lemma_c_t_upper_and_two_lower_bounds}
\end{lemma}
\proof{}
Follows from Lemma~\ref{lemma_c_t_in_terms_of_p_t} and Lemma~\ref{lemma_p_t_upper_and_two_lower_bounds}.
\qed{}

\subsection{Variance of the number of cycles}
\label{proofs_var}

\begin{lemma}
If there exists a function $M = M(n)$ such that, for some $\varepsilon > 0$,
\[
\begin{cases}
M \ge N_1, \\
M p > 1 + \varepsilon, \\
t^2 / M \to 0,
\end{cases}
\]
then
\[\frac{\op{Var}X}{{{(\op{E}X)}^{2}}} \le \frac{1}{\op{E}X} + o\left(\frac{M^t}{t c_{t}}\right).\]
\label{lemma_variance_bound}
\end{lemma}
\proof{}
Let $\gamma_1, \ldots , \gamma_{c_t}$ be all copies of $C_t$ in $G(n, r, s)$ and $X_i$ be the indicator of the event $\left\{\gamma_i \subset G_p(n,r,s)\right\}$. Then
\[\begin{gathered}
\op{Var}X=\sum\limits_{i=1}^{{{c}_{t}}}{\op{Var}{{X}_{i}}}+\sum\limits_{i\ne j}{\op{cov}({{X}_{i}},{{X}_{j}})}=\sum\limits_{i=1}^{{{c}_{t}}}{\op{Var}{{X}_{i}}}+\sum\limits_{(i,j)\in I}{\op{cov}({{X}_{i}},{{X}_{j}})} \le\\
\le \sum\limits_{i=1}^{{{c}_{t}}}{\op{E}(X_{i}^{2})}+\sum\limits_{(i,j)\in I}{\op{E}({{X}_{i}}{{X}_{j}})}=\sum\limits_{i=1}^{{{c}_{t}}}{\op{E}{{X}_{i}}}+\sum\limits_{(i,j)\in I}{\op{E}({{X}_{i}}{{X}_{j}})}=\op{E}X+\sum\limits_{(i,j)\in I}{\op{E}({{X}_{i}}{{X}_{j}})},\end{gathered}\]
where $I$ is the set of ordered pairs $(i, j)$ of indices such that $\gamma_i$ and $\gamma_j$ have at least one common edge.

Let $(i,j) \in I$. Let $\alpha(i,j)$ denote the number of all inclusion maximal non-degenerate (with at least one edge) simple paths in $\gamma_i \cap \gamma_j$. Let $x(i,j)$ denote the total number of edges in $\gamma_i \cap \gamma_j$. Clearly, $\alpha(i,j) \le x(i, j) \le t - \alpha(i, j)$ and $\alpha(i,j) \le t/2$.

Given $\alpha \le \floor{t / 2}$ and $\alpha \le x \le t - \alpha$, let's estimate from above the number of pairs $(i,j) \in I$ such that $\alpha(i,j) = \alpha$ and $x(i,j) = x$. Let's fix the index $i$, which can be done in $c_{t}$ ways, and count the number of ways to choose $j$. For any $j$ such that $\alpha(i,j) = \alpha$ and $x(i,j) = x$, let us write down the vertices of $\gamma_t$ sequentially, following its edges:
\[\begin{aligned}
(& v_{\xi_1}, \ldots, v_{\xi_1 + b_1 \cdot l_1},
u_{1}, \ldots, u_{r_1}, \\
& v_{\xi_2}, \ldots, v_{\xi_2 + b_2 \cdot l_2},
u_{r_1 + 1}, \ldots, u_{r_1 + r_2}, \\
& \ldots, \\
& v_{\xi_{\alpha - 1}}, \ldots, v_{\xi_{\alpha - 1} + b_{\alpha - 1} \cdot l_{\alpha - 1}},
u_{r_1 + \ldots + r_{\alpha - 2} + 1}, \ldots, u_{r_1 + \ldots + r_{\alpha - 1}}, \\
& v_{\xi_{\alpha}}, \ldots, v_{\xi_{\alpha} + b_{\alpha} \cdot l_{\alpha}},
u_{r_1 + \ldots + r_{\alpha - 1} + 1}, \ldots, u_{r_1 + \ldots + r_{\alpha}}).
\end{aligned}\]
Here $v_1, \ldots, v_{t}$ are the vertices of $\gamma_i$ such that $\{v_i, v_{i+1}\} \in E$, $(v_{-(t-1)}, \ldots, v_{0})= (v_1, \ldots, v_{t})$,
and, for $m \in \irange{1}{\alpha}$,
\[\begin{cases}
1 \le \xi_m \le t, \\
b_m \in \{-1, 1\}, \\
1 \le l_m \le l_1 + \ldots + l_{\alpha} = x, \\
0 \le r_m \le r_1 + \ldots + r_{\alpha} = t - (x + \alpha).
\end{cases}\]
Thus, $j$ is uniquely defined by the tuple
\[\left(\xi_1, \ldots, \xi_{\alpha}, b_1, \ldots, b_{\alpha}, l_1, \ldots, l_{\alpha}, r_1, \ldots, r_{\alpha}, u_1, \ldots, u_{t - (x + \alpha)}\right),\]
the number of ways to choose which, considering the constraints, is bounded from above by
\[{{t}^{\alpha }}{{2}^{\alpha }}\binom{x-1}{\alpha -1}\binom{t-x-1}{\alpha -1}N_{1}^{t-x-\alpha}.\]
Note that $N_{1}^{t-x-\alpha}$ is an upper bound for the number of choices of $u_1, \ldots, u_{t - (x + \alpha)}$ if the other parameters are fixed. Indeed, if $u_1, \ldots, u_{t - (x + \alpha)}$ are chosen sequentially, then, for each $k\in\irange{1}{t - (x + \alpha)}$, either $u_k \in V_s(v_{\xi_m + b_m \cdot l_m})$ for some $m\in\irange{1}{\alpha}$ or $u_k \in V_s(u_{k-1})$.

Let's fix some $\delta \in (0,1/2)$ and let $D = M p$. Notice that $D>1+\varepsilon$ is assumed in the conditions in the statement of Lemma~\ref{lemma_variance_bound}. Since the probability of appearance of two fixed cycles with $x$ common edges is $p^{2t - x}$, the sum $\sum\limits_{(i,j)\in I}{\op{E}({{X}_{i}}{{X}_{j}})}$ can be estimated from above as follows:
\[\sum\limits_{(i,j)\in I}{\op{E}({{X}_{i}}{{X}_{j}})}\le 
\sum\limits_{\alpha =1}^{\floor{t/2}} {\sum\limits_{x=\alpha}^{t-\alpha } {c}_{t}{{{(2t)}^{\alpha }}\binom{x-1}{\alpha -1}\binom{t-x-1}{\alpha -1}N_{1}^{t-x-\alpha }{{p}^{2t-x}}}} \le\]
\[\le \frac{1}{t} {c}_{t}{{p}^{t}} \sum\limits_{\alpha=1}^{{\floor{t/2}}} {\sum\limits_{x=\alpha}^{t - \alpha}t{{{(2t)}^{\alpha}}\binom{x}{\alpha}\binom{t-x}{\alpha -1}M^{t-x-\alpha }{{p}^{t-x}}}} =\]
\[= \frac{1}{t} c_t p^t {D^{t}}\sum\limits_{\alpha =1}^{\floor{t/2}}{\sum\limits_{x=\alpha}^{t-\alpha }{{D^{-x}}{{\left( \frac{2t}{{M}} \right)}^{\alpha }}\binom{x}{\alpha}\frac{t \alpha}{t-x-\alpha+1}\binom{t-x}{\alpha }}} \le \]
\[\le \frac{1}{t} c_t p^t D^{t}\sum\limits_{\alpha =1}^{\floor{t/2}}{\sum\limits_{x=\alpha}^{t-\alpha }{{D^{-x}}{{\left( \frac{2t}{{M}} \right)}^{\alpha }}\frac{t \alpha}{t-x-\alpha+1} {{\left( \frac{e x}{\alpha } \right)}^{\alpha }}{{\left(\frac{e(t-x)}{\alpha } \right)}^{\alpha }}}}\le \]
\[\le \frac{1}{t} c_t p^t {{D}^{t}}\sum\limits_{\alpha, x=1}^{+\infty}{{{{D}^{-x}} \frac{t}{\max \left\{1, t-x-\alpha+1 \right\}} {{\left(\alpha^{1 / \alpha} \frac{2{{e}^{2}}{{t}^{2}}}{M} \frac{x}{{{\alpha }^{2}}} \right)}^{\alpha }}}}\le\]
\[\le \frac{1}{t} c_t p^t {{D}^{t}}{\sum\limits_{\alpha, x=1}^{+\infty}{\frac{t (1 + \varepsilon)^{-\delta (x + \alpha)}}{\max \left\{1, t-x-\alpha+1 \right\}} {{D}^{-x (1 - \delta)}} {{\left((1 + \varepsilon)^{\delta} \alpha^{1 / \alpha} \frac{2{{e}^{2}}{{t}^{2}}}{M} \frac{x}{{{\alpha}^{2}}} \right)}^{\alpha}}}}.\]
Further, for $x + \alpha \le \frac{t}{2}$,
\[\frac{t (1 + \varepsilon)^{-\delta (x + \alpha)}}{\max \left\{1, t-x-\alpha+1 \right\}} \le \frac{1}{1-\frac{x+\alpha-1}{t}} \le \frac{1}{1-1/2} = 2\]
and, for $x + \alpha > \frac{t}{2}$,
\[\frac{t (1 + \varepsilon)^{-\delta (x + \alpha)}}{\max \left\{1, t-x-\alpha+1 \right\}} \le t (1 + \varepsilon)^{-\delta t / 2} \le \frac{2}{\delta\ln (1 + \varepsilon)}.\]
Thus, in any case,
\[\frac{t (1 + \varepsilon)^{-\delta (x + \alpha)}}{\max \left\{1, t-x-\alpha+1 \right\}} \le 2 + \frac{2}{\delta\ln (1 + \varepsilon)}.\]
Let
\[d = (1 + \varepsilon)^{\delta} \frac{2{{e}^{2}}{{t}^{2}}}{M} \underset{\alpha}{\max}\left(\alpha\left(2 + \frac{2}{\delta\ln (1 + \varepsilon)}\right)\right)^{1/\alpha}.\]
Then
\[\sum\limits_{(i,j)\in I}{\op{E}({{X}_{i}}{{X}_{j}})} \le \frac{1}{t} c_t p^t {{D}^{t}}{\sum\limits_{\alpha, x=1}^{+\infty}{{{D}^{-x (1 - \delta)}} {{\left(d \frac{x}{{{\alpha }^{2}}} \right)}^{\alpha }}}} \le \frac{1}{t} c_t p^t {{D}^{t}}{\sum\limits_{\alpha, x=1}^{+\infty}{(1 + \varepsilon)^{-\delta x} f(x)}},\]
where \[f(x) = {{{D}^{-x (1 - 2\delta)}} {{\left(d \frac{x}{{{\alpha }^{2}}} \right)}^{\alpha }}}.\]
From the equation $\partial \ln f / \partial x = 0$, $f$ achieves its maximum at $x_{max} = \frac{\alpha}{(1 - 2\delta)\ln D}$. Since, under Lemma~\ref{lemma_variance_bound} assumptions, $d\to 0$ as $n\to\infty$, then, for sufficiently large $n$,
\[f(x_{max}) = \left(\frac{d}{e \alpha (1 - 2\delta) \ln D}\right)^{\alpha} \le d^{\alpha/2}.\]
Finally, for large enough $n$,
\[\begin{gathered}
\sum\limits_{(i,j)\in I}{\op{E}({{X}_{i}}{{X}_{j}})} \le \frac{1}{t} c_t p^t {{D}^{t}}\sum\limits_{\alpha, x=1}^{+\infty} (1 + \varepsilon)^{-\delta x} d^{\alpha/2} = \\
= \frac{1}{t} c_t p^t D^t \cdot \frac{1}{(1 + \varepsilon)^{\delta}-1} \cdot \frac{d^{1/2}}{1 - d^{1/2}} = o \left( \frac{1}{t} c_t p^t D^t \right)
\end{gathered}\]
and
\[\frac{\op{Var}X}{{{(\op{E}X)}^{2}}} \le \frac{\op{E}X + o \left( \frac{1}{t} c_t p^t D^t \right)}{{{(\op{E}X)}^{2}}} = \frac{1}{\op{E}X} + o\left(\frac{M^t}{t c_{t}}\right).\]
\qed{}

\subsection{Proof of Theorem~\ref{th_long_thresh}}
\label{proofs_proof}
As before, $X$ is the number of copies of $C_t$ in $G_p(n, r, s)$.
Let $p = o\left(\hat{p}\right)$. Then by Markov's inequality and Lemma~\ref{lemma_c_t_upper_and_two_lower_bounds}:
\[\mathbb{P}(X \ge 1) \le \op{E}X = c_t p^t = O \left(\frac{1}{2t} (N_1pn^{s/t})^{t}\right) = O \left(\frac{1}{2t} (p / \hat{p})^{t}\right) \to 0.
\]
\bigskip
Let $p = \omega\left(\hat{p}\right)$. Then Lemma~\ref{lemma_c_t_upper_and_two_lower_bounds} and Lemma~\ref{lemma_variance_bound} for $M = N_1 \cdot \max \left\{1, \left.n^{s/t} \middle/ \binom{r}{s}\right.\right\}$ yield
\[\op{E}X = \Omega\left(\frac{1}{2t} \left(N_1 p n^{s/t} \middle/ \binom{r}{s}\right)^{t}\right) = \frac{1}{2t} \left(\omega(1)\right)^{t} \to +\infty,\]
\[\frac{\op{Var}X}{{{(\op{E}X)}^{2}}} \le \frac{1}{\op{E}X} + o\left(\frac{N_1^t \cdot \max \left\{1, \left.n^{s} \middle/ \binom{r}{s}^t\right.\right\}}{t \cdot \max\left\{\frac{1}{2t} N_1^{t}, \frac{n^s}{2t} \left(N_1 \middle/ \binom{r}{s}\right)^{t}\right\}} \right) \to 0.\]
By Chebyshev's inequality:
\[\mathbb{P}(X = 0) \le \mathbb{P}\left(|X - \op{E}X| \ge \op{E}X \right) \le \frac{\op{Var}X}{{{(\op{E}X)}^{2}}} \to 0.\]
Thus, it is proved that $\hat{p}$ is a threshold.
\bigskip

Now, it remains to prove that the threshold is sharp if $s=0$ and $t\to+\infty$ or $s$ is arbitrary and $t=\omega(\ln n)$.
If $t \to +\infty$, then, for $p \le (1 - \varepsilon) \hat{p}$, we get
\[\op{E}X = c_t p^t = O \left(\frac{1}{2t} (N_1pn^{s/t})^{t}\right) = O \left(\frac{1}{2t} (1-\varepsilon)^{t}\right) \to 0.
\]
If, in addition, $s=0$ or $t = \omega(\ln n)$, then for $p \ge (1 + \varepsilon) \hat{p}$,
\[N_1 p \ge (1 + \varepsilon) n^{-s/t} \sim 1 + \varepsilon > 1 + \frac{\varepsilon}{2}.\]
Therefore, Lemma~\ref{lemma_c_t_upper_and_two_lower_bounds} and Lemma~\ref{lemma_variance_bound} for $M = N_1$ yield
\[\op{E}X = c_t p^t = \Omega\left(\frac{1}{2t} (N_1p)^{t}\right) = \Omega\left(\frac{1}{2t} \left(1+\frac{\varepsilon}{2}\right)^t\right) \to +\infty,
\]
\[\frac{\op{Var}X}{{{(\op{E}X)}^{2}}} \le \frac{1}{\op{E}X} + o\left(\frac{N_1^t}{t c_{t}}\right) = \frac{1}{\op{E}X} + o(1) \to 0\]
which means that the threshold is sharp.
\qed{}

\section{Conclusion}

We have established the threshold for containment of a copy of $C_t$ in $G_p(n,r,s)$ under assumption that $t=o\left(\sqrt{N_1}\right)$. This restriction is due to the method we used to prove auxiliary results about the number of copies of $C_t$ in $G(n,r,s)$ (Lemma~\ref{lemma_c_t_upper_and_two_lower_bounds}) and the variance of the number of copies of $C_t$ in $G_p(n,r,s)$ (Lemma~\ref{lemma_variance_bound}).

In fact, the probability in \eqref{eq_sharp_thresh} is still a threshold for $t=O\left(\sqrt{N_1}\right)$. It can be shown using the same techniques and a bit more careful analysis. The bounds for $p_t$ and $c_t$ in Lemmas~\ref{lemma_p_t_upper_and_two_lower_bounds} and \ref{lemma_c_t_upper_and_two_lower_bounds} still hold for $t=O\left(\sqrt{N_1}\right)$. The conditions in Lemma~\ref{lemma_variance_bound} would be
\[\begin{cases}
	M \ge N_1, \\
	M p \to \infty, \\
	t^2 / M = O(1).
\end{cases}\]
Note that if $t^2 / M = \Theta(1)$, the condition $M p > 1 + \varepsilon$ is insufficient to prove Lemma~\ref{lemma_variance_bound} using the same idea. This means that we do not know whether the threshold~\eqref{eq_sharp_thresh} is sharp for $t=\Theta\left(\sqrt{N_1}\right)$.

\section{Acknowledgements}

The reported study was funded by RFBR, project number 20-24-70001, and by Grant N NSh-2540.2020.1 to support leading scientific schools of Russia.


\begin{thebibliography}{99}

\bibitem{Bollobas1981} B.~Bollob\'{a}s. Threshold functions for small subgraphs. Mathematical Proceedings of the Cambridge Philosophical Society, 1981, vol.~90, no.~2, pp.~197--206.

\bibitem{Bollobas1987} B.~Bollob{\'{a}}s, A.G.~Thomason. Threshold functions. Combinatorica, 1987, vol.~7, no.~1, pp.~35--38.

\bibitem{Bollobas1989} B.~Bollob\'{a}s, J.~Wierman. Subgraph Counts and Containment Probabilities of Balanced and Unbalanced Subgraphs in a Large Random Graph. Annals of the New York Academy of Sciences, 1989, vol.~576, no.~1, pp.~63--70

\bibitem{Bollobas2016} B.~Bollob\'{a}s, B.P.~Narayanan, A.M.~Raigorodskii. On the stability of the Erd{\H{o}}s-Ko-Rado theorem. Journal of Combinatorial Theory, Series A, 2016, vol.~137, pp.~64--78.

\bibitem{Burkin2016} A.V.~Burkin. Small Subgraphs in Random Distance Graphs. Theory of Probability {\&} Its Applications, 2016, vol.~60, no.~3, pp.~367--382.

\bibitem{Burkin2018} A.V.~Burkin, M.E.~Zhukovskii. Small subgraphs and their extensions in a random distance graph. Sbornik: Mathematics, 2018, vol.~209, no.~2, pp.~163--186.

\bibitem{Erdos1959}
P.~Erd\H{o}s, A.~R\'{e}nyi. On random graphs. I. Publicationes Mathematicae Debrecen, 1959, vol.~6, pp.~290--297.

\bibitem{Erdos1960} P.~Erd\H{o}s, A.~R\'{e}nyi. On the evolution of random graphs. Publications of the Mathematical Institute of the Hungarian Academy of Sciences, 1960, vol.~5, pp.~17--61.

\bibitem{Friedgut2006} E.~Friedgut, V.~R\"{o}dl, A.~Ruci{\'{n}}ski, P.~Tetali. A sharp threshold for random graphs with a monochromatic triangle in every edge coloring. Memoirs of the American Mathematical Society, 2006, vol.~179, no.~845.

\bibitem{Gilbert1959} E.N.~Gilbert. Random Graphs. The Annals of Mathematical Statistics, 1959, vol.~30, no.~4, pp.~1141--1144.

\bibitem{Janson2000} S.~Janson, T.~Luczak, A.~Ruci{\'{n}}ski. Random Graphs, 2000, Wiley, New York.

\bibitem{Karoski1983} M.~Karo{\'{n}}ski, A.~Ruci{\'{n}}ski. On the number of strictly balanced subgraphs of a random graph. Graph Theory. Lecture Notes in Mathematics, 1983, vol.~1018, pp.~79--83

\bibitem{Komlos1983} J.~Koml{\'{o}}s, E.~Szemer{\'{e}}di. Limit distribution for the existence of hamiltonian cycles in a random graph. Discrete Mathematics, 1983, vol.~43, no.~1, pp.~55--63.

\bibitem{Kupavskii2016} A.~Kupavskii. On random subgraphs of Kneser and Schrijver graphs. Journal of Combinatorial Theory, Series A, 2016, vol.~141, pp.~8--15.

\bibitem{Kupavskii2018} A.~Kupavskii. Random Kneser Graphs and Hypergraphs. The Electronic Journal of Combinatorics, 2018, vol.~25, no.~4.

\bibitem{Ogarok2020} P.A.~Ogarok, A.M.~Raigorodskii. On Stability of the Independence Number of a Certain Distance Graph. Information Transmission Problems, 2020, vol.~56, no.~4, pp.~347--359.

\bibitem{Pyaderkin2016} M.M.~Pyaderkin. Independence numbers of random subgraphs of distance graphs. Mathematical Notes, 2016, vol.~99, no.~3-4, pp.~556--563.

\bibitem{Raigorodskii2021} A.M.~Raigorodskii, V.S.~Karas. Asymptotics of the independence number of a random subgraph of the graph $G(n,r,<s)$. Submitted to Mathematical Notes.

\bibitem{Rucinski1986} A.~Ruci{\'{n}}ski, A.~Vince. Strongly balanced graphs and random graphs. Journal of Graph Theory, 1986, vol.~10, no.~2, pp.~251--264.

\bibitem{Rucinski1988} A.~Ruci{\'{n}}ski. When are small subgraphs of a random graph normally distributed? Probability Theory and Related Fields, 1988, vol.~78, no.~1, pp.~1--10.

\bibitem{Schurger1976} K.~Sch\"{u}rger. On the evolution of random graphs over expanding square lattices. Acta Mathematica Academiae Scientiarum Hungaricae, 1976, vol.~27, no.~3-4, pp.~281--292.

\bibitem{Zhukovskii2010} M.E.~Zhukovskii. The weak zero-one law for random distance graphs. Teoriya Veroyatnostei i ee Primeneniya (in Russian), 2010, vol.~55, no.~2, pp.~344--350. [Theory of Probability and Its Applications (English translation), 2010, vol.~55, no.~2, pp.~356--360].

\bibitem{Zhukovskii2012_zero} M.E.~Zhukovskii. A weak zero-one law for sequences of random distance graphs. Sb. Math., 2012, vol.~203, no.~7, pp.~1012--1044.

\bibitem{Zhukovskii2012_sub} M.E.~Zhukovskii. On the probability of the occurrence of a copy of a fixed graph in a random distance graph. Mathematical Notes, 2012, vol.~92, no.~5-6, pp.~756--766.


\end{thebibliography}
\end{document}